\documentclass[12pt]{article}
\usepackage{latexsym,amsmath}

\newcommand{\qed}{\hfill $\Box$ \\}
\def\deg{{\rm deg}}
\def\e{\epsilon}

\def\cF{{\cal F}}

\def\cI{{\cal I}}

\def\cH{{\cal H}}

\newcommand{\brac}[1]{\left(#1\right)}
\newcommand{\bfrac}[2]{\brac{\frac{#1}{#2}}}
\newcommand{\beq}[1]{\begin{equation}\label{#1}}
\newcommand{\eeq}{\end{equation}}
\newcommand{\blem}[1]{\begin{lemma}\label{#1}}
\newcommand{\elem}{\end{lemma}}
\newcommand{\bth}[1]{\begin{theorem}\label{#1}}
\newcommand{\enth}{\end{theorem}}
\newcommand{\brem}[1]{\begin{remark}\label{#1}}
\newcommand{\erem}{\end{remark}}

\def\cF{{\cal F}}

\def\a{\alpha}

\def\e{\varepsilon}

\def\th{\theta}

\def\m{\mu}

\def\r{\rho}

\def\cF{{\cal F}}

\def\whp{{\bf whp}}

\newcommand{\set}[1]{\left\{#1\right\}}

\newcommand{\proofend}{\hspace*{\fill}\mbox{$\Box$}}

\def\E{{\sf E}}
\def\Pr{{\sf P}}
\def\cE{{\cal E}}

\begin{document}
\newtheorem{theorem}{Theorem}
\newtheorem{corollary}[theorem]{Corollary}
\newtheorem{lemma}[theorem]{Lemma}
\newtheorem{proposition}[theorem]{Proposition}
\newtheorem{conjecture}[theorem]{Conjecture}
\newtheorem{definition}[theorem]{Definition}
\newtheorem{claim}{Claim}
\newtheorem{remark}{Remark}

\title{Coloring $H$-free hypergraphs}
\author{Tom Bohman\thanks{Department of Mathematical Sciences,
Carnegie Mellon University,  Pittsburgh, PA 15213-3890; Supported in
part by  NSF grant DMS 0701183},\, Alan Frieze\thanks{Department of
Mathematical Sciences, Carnegie Mellon University,  Pittsburgh, PA
15213-3890; Supported in part by  NSF grant CCF0502793}, \,Dhruv
Mubayi\thanks{ Department of Mathematics, Statistics, and Computer
Science, University of Illinois, Chicago, IL 60607-7045; Supported
in part by  NSF grant DMS 0653946} }

\maketitle

\vspace{-0.4in}

\begin{abstract}
Fix $r \ge 2$ and a collection of $r$-uniform hypergraphs $\cH$. What is the minimum number of edges in an $\cH$-free $r$-uniform hypergraph
with chromatic number greater than $k$? We investigate this question for various $\cH$. Our results include the following:

\noindent
$\bullet$ An $(r,l)$-system is an $r$-uniform hypergraph with every two edges sharing at most $l$ vertices. For $k$ sufficiently large, there is
an $(r,l)$-system
with chromatic number greater than $k$ and number of edges at most $c(k^{r-1}\log k)^{l/(l-1)}$, where
$$c=2\bfrac{100(r)_l^2}{l!}^{1/(l-1)}\bfrac{10(r-1)}{l-1}^{l/(l-1)}.$$
This improves on the previous best bounds of Kostochka-Mubayi-R\"odl-Tetali~\cite{KMRT}. The upper bound is sharp apart from the constant $c$ as
shown in~\cite{KMRT}.

\noindent
$\bullet$
The minimum number of edges in an $r$-uniform hypergraph with independent neighborhoods and chromatic number greater than $k$ is of order
$k^{r+1/(r-1)}\log^{O(1)}k$ as $k \rightarrow \infty$. This generalizes (aside from logarithmic factors) a result of Gimbel and Thomassen~\cite{GT} for
triangle-free graphs.

\noindent
$\bullet$
Let $T$ be an $r$-uniform  hypertree of $t$ edges. Then every $T$-free $r$-uniform hypergraph has chromatic number at most $2(r-1)(t-1)+1$.
This generalizes the well known fact that every $T$-free graph has chromatic number at most $t$.

Several open problems and conjectures are also posed.

\end{abstract}

\section{Introduction}

An $r$-graph is a hypergraph whose edges all have size $r$. The chromatic number of an $r$-graph is the minimum number of colors required to partition its vertex set so that no edge is monochromatic. The starting point of our investigations is the following basic question:
\smallskip

What is the minimum number $m_k(r)$ of edges in an $r$-graph with chromatic number greater than $k$?
\smallskip

In general this problem is very difficult to solve exactly, and so we seek asymptotic results as one or both of $k,r$ tend to infinity.  It is easy to observe that $m_k(2)=\binom{k+1}{2}$, but already determining $m_k(3)$ is a challenging open question: $m_2(3)=7$ achieved by the Fano plane, but $m_3(3)$ is unknown.
For fixed $r$ and large $k$, the best known bounds are due to Alon and are still far apart.

\begin{theorem} {\bf (Alon \cite{A})}\label{alon}
For every $k,r \ge 2$,
$$(r-1)\left\lceil\frac{k}{r}\right\rceil\left(\frac{r-1}{r} k\right)^{r-1}<m_k(r)
<\binom{kr+1}{r}\frac{5\log r}{r}.$$
\end{theorem}
Note that this implies that for fixed $r$ and $k \rightarrow \infty$ we have $m_k(r)=\Theta(k^r)$. In the opposite direction, determining $m_2(r)$
is the well known problem concerning the minimum number of edges in $r$-graphs of chromatic number greater than two, generally referred
to as not having Property B.
The best upper bound follows from an old probabilistic construction of Erd\H os~\cite{EpropB} while Radhakrishnan and Srinivasan~\cite{RS} proved the lower bound below (for large $r$) which is the best to date.
$$0.7 \sqrt {r/\log r}\times 2^r<m_2(r)  < r^22^r.$$

Here we consider the same question, but we impose a natural restriction on the underlying $r$-graph.

\begin{definition}
Fix $r \ge 2$ and a collection of $r$-graphs $\cH$. Let $m_k(\cH)$ be the minimum number of edges in an $\cH$-free $r$-graph with chromatic number greater than $k$.
\end{definition}
Note that $\cH$-free in the definition above refers to not necessarily induced subhypergraphs. Also, if $\cH=\{H\}$, we will abuse notation by writing $m_k(H)$.

Our goal is to determine $m_k(\cH)$ for various $\cH$. Special cases
of this parameter have already been studied, and lead to difficult
problems.  For example, Gimbel and Thomassen~\cite{GT} proved that
$m_k(K_3)$ has order of magnitude $k^3\log^2 k$ as $k \rightarrow
\infty$. However, determining $m_k(K_4)$ and $m_k(C_4)$ are open
problems. In fact, these problems seem harder than determining the
Ramsey numbers for the corresponding graphs, and the growth rates of
these Ramsey numbers are well-studied and not known. Call a
hypergraph {\em nontrivial} if it has at least two edges.

\begin{definition}
Let $H$ be a nontrivial $r$-graph. Let
$$\r(H)=\max_{H' \subset H}\,\frac{e'-1}{v'-r},$$
where $H'$ is nontrivial with $v'$ vertices and $e'$ edges.
For a finite family $\cH$ of nontrivial $r$-graphs,
$\r(\cH)=\min_{H \in \cH} \r(H)$.
\end{definition}
We say that $H$ is {\em balanced} if this maximum occurs for $H$ i.e.
$$\r(H)=\frac{e(H)-1}{v(H)-r}.$$

The parameter $\r$ appears to be the crucial hypergraph invariant for our problem. Our main result stated below provides a very general upper
bound for $m_k(\cH)$. As we will show, in many cases this general upper bound seems to give the correct order of magnitude for fixed $r$ as $k
\rightarrow \infty$.

\begin{theorem} \label{main}
Let $\cH=\{H_1,H_2,\ldots,H_\ell\}$ be a finite family of nontrivial
balanced $r$-graphs with $v_i=v(H_i)$ and $e_i=e(H_i)$. Let
$\r_i=\r(H_i)$ and assume that $\r_i\leq \r_{i+1}$ for $1\leq i<
\ell$. Define $s$ by $\r=\r_1=\r_2=\cdots=\r_s<\r_{s+1}$ and assume
that $\r>1/(r-1)$. For each $i$ and each edge $e\in H_i$ let
$\a_i(e)$ be the number of automorphisms of $H_i$ that map $e$ to
itself.
Let $\a_i=\min_e\a_i(e)$.\\
Suppose that $c_1$ is the solution to
$$\sum_{i=1}^s\frac{e_ix^{e_i-1}}{\a_i}=\frac{1}{50r!}.$$
and let
$$c_2=\max\set{\bfrac{\a_i}{5e_ic_1^{e_i-1}r!}^{1/(r-1)}:\;i=1,2,\ldots,s}\ge 10^{1/(r-1)}.$$
Then for large $k$,
$$m_k(\cH)<c_{\cH} (k^{r-1}\log k)^{(r-1/\r)/(r-1-1/\r)}$$
where
$$c_{\cH}=2\bfrac{r!}{c_1}^{1/(r-1-{1/\r})}c_2^{(r-1)(r-{1/\r})/(r-1-{1/\r}))}
\bfrac{r-1}{r-1-{1/\r}}^{(r-{1/\r})/(r-1-{1/\r})}$$
(We can if we wish replace the 2 by a constant arbitrarily close to 1).
\end{theorem}
Note that the exponent of $k$ in Theorem \ref{main} is always greater than $r$.

\begin{remark}
The restriction to a collection of balanced hypergraphs is not too
restrictive. Our applications will be balanced and in general, if we
replace an $H$ by a subhypergraph $H'$ that determines $\r(H)$ then
the upper bound we obtain for $m_k(\cH)$ is valid. After all, a
hypergraph that does not contain $H'$ cannot contain $H$.
\end{remark}

\subsection{$(r,l)$-systems}
An $(r,l)$-system is an $r$-graph with every pair of edges sharing fewer than $l$ vertices. Let $m_k(r,l)$ denote the minimum number of edges in
an $(r,l)$-system with chromatic number greater than $k$. Erd\H os and Lov\'asz~\cite{EL} studied $m_k(r,2)$, indeed the Local lemma was
originally developed and used to give lower bounds for this parameter.  Recently, Kostochka et.al.~\cite{KMRT} proved that
$m_k(r,l)$ has order of magnitude $(k^{r-1}\log k)^{l/(l-1)}$ as $k \rightarrow \infty$. They proved the upper bound
$m_k(r,l)<b_{r,l}(k^{r-1}\log k)^{l/(l-1)}$ where
$$b_{r,l}=\frac{2(2r^{3l})^{l/(l-1)}}{(r)_l}.$$
 Using Theorem \ref{main} we can substantially improve this constant.

\begin{theorem} \label{rl}
Fix $2\le l <r$ and let $k$ be sufficiently large. Then $m_k(r,l)<c_{r,l}(k^{r-1}\log k)^{l/(l-1)}$, where
$$c_{r,l}=2\bfrac{100(r)_l^2}{l!}^{1/(l-1)}\bfrac{10(r-1)}{l-1}^{l/(l-1)}.$$
\end{theorem}
Note that for large $r$, $b_{r,l}$ grows like $r^{2l}$ whereas $c_{r,l}$ grows like $r^3$.
\subsection{Independent neighborhoods}
A triangle-free graph is one whose neighborhoods are all independent sets. Generalizing to $r$-graphs, one can study $r$-graphs with independent
neighborhoods.
If $S$ is a set of vertices in an $r$-graph $G=(V,E)$ and $|S|=r-1$, then its neighborhood
$N_G(S)=\set{v\in V-S: S \cup \{v\} \in E}$.
The degree $\deg_G(S)=|N_G(S)|$.

An $r$-graph has independent neighborhoods if it contains no copy of $F_r$, where $F_r$ is the
$r$-graph comprising $r+1$ edges $\set{E_0,E_1,\ldots,E_r}$. Here, if $A=\cap_{i=1}^r E_i$ then (i) $|A|=r-1$ and
(ii) $E_0=\bigcup_{i=1}^r(E_i\setminus A)$. Thus $F_2=K_3$. Gimbel and Thomassen~\cite{GT} proved that the order of magnitude of $m_k(F_2)$ is
$k^3(\log k)^2$. Although we are unable to determine the correct logarithmic factors, we generalize this result as follows; the upper bound
follows directly from Theorem \ref{main}.

\begin{theorem} \label{indnbd}
Fix $r \ge 3$ and let $k$ be sufficiently large. The minimum number of edges in an $r$-graph with independent neighborhoods and chromatic number
greater than $k$ satisfies
$$b_\cI k^{r+1/(r-1)}<m_k(F_r)<c_\cI k^{r+1/(r-1)}(\log k)^{1+r/(r-1)^2},$$
where
\begin{eqnarray*}
b_\cI&=&\frac{1}{40r^2 2^r}\\
c_\cI&=&\brac{r!\bfrac{50r!(r+1)}{(r-1)!^2}^{1/r}}^{1/(r-3/2)}\bfrac{10(r-1)}{r-3/2}^{(r-1/2)/(r-3/2)}.
\end{eqnarray*}
\end{theorem}
Note that as $r\to\infty$, $c_\cI=O(r)$.

It is hard to even make a conjecture about the correct growth rate of $m_k(F_r)$. Most likely neither
the upper nor lower bounds give the correct order of magnitude. However, improving either bound seems difficult, since the corresponding
improvement for the graph case involved deep results of Kim \cite{Kim} and Johansson \cite{Joh} on the independence number and chromatic number
of triangle-free graphs.  Currently, hypergraph versions of these two results do not exist.

\subsection{Excluding a hypertree}
 A cycle of length $t\ge 2$ in an $r$-graph is a collection of $t$ distinct vertices
 $X=\{x_1, \ldots, x_t\}$ and $t$ distinct edges $E_1,
\ldots, E_t$ such that $\set{x_i,x_{i+1}} \subset E_i$ for each
$i=1, \ldots, t$ (indices taken modulo $t$). An $r$-forest is an
$r$-graph with no cycles. It is easy to see that if $H$ contains a
cycle, then $\rho(H)>1/(r-1)$ and Theorem~\ref{main} applies.  On
the other hand, if $H$ is an $r$-forest, then it is easy to show
that every $H$-free $r$-graph $G$ has chromatic number at most
$c_H$, so there can be no analogue of the upper bound in Theorem
\ref{main}. It is an easy exercise to produce a proper coloring of
$G$ where the number of colors is exponential in the size of $H$.
  The next Theorem shows that we can reduce this bound substantially. An $r$-tree is a connected $r$-forest, where connected means that for
every two vertices $x,y$, there is a sequence of edges $E_1, \ldots E_l$ such that $x \in E_1$, $y \in E_l$, and $E_i \cap E_{i+1}\ne \emptyset$
for all $i=1,2, \ldots, l-1$.  The statement below applies to $r$-trees,
but a similar statement can be proved for $r$-forests as well.

\begin{theorem} \label{Tmain}  Let $T$ be an $r$-tree with $(r-1)t+1$ vertices  and suppose that $G$ is an $r$-graph not containing $T$.  Then
the chromatic number of $G$ is at most $2(r-1)(t-1)+1$.
\end{theorem}

When $r=2$ it is a well-known fact that every $T$-free graph $G$ has
chromatic number at most $t$, when $T$ has $t$ edges. Indeed, this
follows from the observation that every subgraph of $G$ has a vertex
of degree less than $t$.  For $r \ge 3$ such a statement is false.
For example, let $T$ be the 3-tree comprising three edges, not all
containing the same vertex. Let $G$ be the 3-graph on $n$ vertices,
$n$ large, all of whose edges contain a fixed vertex.   Then clearly
$T \not\subset G$ and $G$ has minimum degree $n-2$, which can be
arbitrarily large.  This is the reason that Theorem \ref{Tmain} is
not trivial. Nevertheless, the best lower bound on the chromatic
number of a $T$-free $r$-graph that we have is $t$. It would be very
interesting to narrow the gap for this problem, and we believe that
Theorem \ref{Tmain} is far from the truth\footnote{Recently Po-Shen
Loh has proved optimal results for this problem}.

\subsection{Graphs vs hypergraphs}
Our final result shows the limitations of Theorem~\ref{main} in the case $r> 2$. Let $K_t^r$ be the complete $r$-graph on $t$ vertices. Then
Theorem~\ref{main} implies that $m_k(K_t^r)<k^{r+\e}$ for some positive $\e$ depending on $r$ and $t$. However, for $r \ge 3$, this can be
improved.

\begin{theorem} \label{cliques}
Fix  $t>r \ge 3$.  Then $m_k(K_{t}^r)=k^{r+o(1)}$, where $o(1)\rightarrow 0$ as $k \rightarrow \infty$. On the other hand, for each $s \ge 3$,
there exists $\e=\e_s>0$ such that $m_k(K_s)>k^{2+\e}$.
\end{theorem}

Theorem \ref{cliques} shows an interesting difference between graphs and hypergraphs. In fact, we
conjecture that  a similar result holds if we forbid much less than a clique.  Call an $r$-graph simple
if every two of its edges share at most one vertex; in the notation of Section 1.1, an $r$-graph is simple if and only if it is an
$(r,2)$-system.  Simple hypergraphs are often studied due their similarity to graphs. We believe that  there exist simple $r$-graphs $H$ such
that $m_k(H)=k^{r+o(1)}$.
For $r \ge 4$ this follows from recent unpublished results of R\"odl-Schacht and the third author, however, this remains open for $r=3$.

\begin{conjecture} \label{simple} There exists a simple 3-graph $H$ for which $m_k(H)=k^{3+o(1)}$.
\end{conjecture}

Let $F$ be the Fano plane, which is the 3-graph with seven vertices and seven edges obtained from the points and lines of the projective
geometry of dimension two over the finite field of order two. Perhaps one can even strengthen Conjecture \ref{simple}  by proving that
$m_k(F)=k^{3+o(1)}$?

In the next section we present the proof of Theorem \ref{main}.  Sections 3, 4,  5 and 6 contain the proofs of Theorems \ref{rl}, \ref{indnbd},
\ref{Tmain} and \ref{cliques} respectively.  The last section has several concluding remarks and open problems.

\section{General upper bound: Proof of Theorem \ref{main}}\label{general}
In this section we prove Theorem \ref{main}.  Our proof uses the method developed by Krivelevich \cite{K}
to obtain bounds for off diagonal Ramsey numbers.  The main idea is to take a random hypergraph with appropriate edge probability and
judiciously delete all copies of $H$ from it.
The additional requirement for us is to keep track of the total number of edges.
\bigskip

\noindent
{\bf {Proof of Theorem \ref{main}.}} Let
$$p=c_1n^{-1/\r}$$
and let
$G_p$ be the random $r$-graph on $n$ vertices with edge probability $p$. Let
$E_p=|E(G_p)|$. Then
\beq{1}
E_p\leq \frac{2c_1n^{r-{1/\r}}}{r!}\ \whp.
\eeq

Next let
$$t=c_2\bfrac{r!\log n}{p}^{1/(r-1)}= c_2\brac{\frac{r!\log n}{c_1}n^{{1/\r}}}^{1/(r-1)}.$$
Now, using the Chernoff bounds to get the first inequality below, we have
\begin{eqnarray*}
\Pr\brac{\exists S:|S|=t\ and\ |E(S)|\leq
  E_0=\binom{t}{r}p/2}&\leq&\binom{n}{t}\exp\set{-\frac{1}{8}\binom{t}{r}p}\\
&\leq&\brac{\frac{ne}{t}\exp\set{-\frac{t^{r-1}}{10r!}p}}^t\\
&=&\brac{\frac{ne}{t}\exp\set{-\frac{1}{10}c_2^{r-1}\log n}}^t\\
&=&o(1).
\end{eqnarray*}
So, \whp:
\beq{2}
\mbox{Every $t$-set contains at least $E_0$ edges.}
\eeq
Now, for $|S|=t$ let
$Y_{S,i}$ be the number of edges in copies of $H_i$ containing at least one
edge of $S$.
Let $Z_{S,i}$ be the number of edges in a maximal collection of pair-wise disjoint copies of $H_i$,
each containing at least one edge of $S$.

Clearly,
$$Z_{S,i}\leq Y_{S,i}.$$
Let
$$\m_i=\binom{t}{r}\binom{n}{v_i-r}e_ip^{e_i}\frac{r!(v_i-r)!}{\a_i}.$$

Thus
$$\E(Y_{S,i})\leq \m_i.$$
{\bf Explanation:}
We choose an edge $e$ of $H_i$ and an $r$-subset $R$ of $S$ to fix an edge that will be $e$ in a copy of
$H_i$. Then we choose $v_i-r$ other vertices for the remainder of our copy. This accounts for $e_i\binom{t}{r}\binom{n}{v_i-r}$
choices. We then choose a copy of $H_i$ in these $v_i$ vertices for which $R$ is an edge. The number of ways of doing this is
$\frac{r!(v_i-r)!}{\a_i(e)}\leq \frac{r!(v_i-r)!}{\a_i}$. Finally, we multiply by $p^{e_i}$, the probability that the $e_i$ edges chosen actually exist.

Now for $A>0$,
$$\Pr(Z_{S,i}\geq A\m_i)\leq
\frac{\E(Z_{S,i})^{A\m_i}}{(A\m_i)!}\leq
\frac{\m_i^{A\m_i}}{(A\m_i)!}\leq \bfrac{e}{A}^{A\m_i}.
$$
(Here we are using an inequality of Erd\H{o}s and Tetali \cite{ET}, see for example Lemma 8.4.1 of \cite{AS}).

Suppose now that
$$A_i=\begin{cases}9&i=1\\10&2\leq i\leq s\\\frac{1}{2^{i-s}}&i>s\end{cases}.$$
Then
\beq{3}
\Pr(\exists S:Z_{S,i}\geq
A_i\m_i)\leq\binom{n}{t}\bfrac{e}{A_i}^{A_i\m_i}.
\eeq
Since $n$ is sufficiently large, we have for some $\frac{9}{10}\leq \th_i\leq 1$,
\begin{eqnarray}
\m_i&=&\frac{\th_ie_ic_1^{e_i}t^rn^{v_i-r-e_i/\r}}{\a_i}\nonumber\\
&=&\frac{\th_ie_ic_1^{e_i-1}c_2^{r-1}r!n^{v_i-r-(e_i-1)/\r}\log n}{\a_i}t\nonumber\\
&=&\frac{\th_ie_ic_1^{e_i-1}c_2^{r-1}r!n^{(v_i-r)(1-\r_i/\r)}\log n}{\a_i}t\label{4}
\end{eqnarray}
If $i\leq s$ then from the definition of $c_2$ we see that $\m_i\geq \frac{1}{5}\th_it\log n$. Hence,
\beq{5}
\Pr(\exists S:Z_{S,i}\geq A_i\m_i)\leq
\brac{\frac{ne}{t}\bfrac{e}{9}^{9\log n/5}}^t=o(1).
\eeq
Now for $i>s$, and because $\r_1=\r$,
$$\frac{A_i\m_1}{\m_i}\geq n^{(v_i-r)(\r_i/\r-1)-o(1)}.$$

So,
\beq{6}
\Pr(\exists S:Z_{S,i}\geq A_i\m_1)\leq \binom{n}{t}\bfrac{e\m_i}{A_i\m_1}^{A_i\m_1}\leq
\binom{n}{t}n^{-\Omega(t\log n)}=o(1).
\eeq

It follows from \eqref{5}, \eqref{6} that \whp
\beq{7}
\sum_{i=1}^\ell Z_{S,i}\leq 10\sum_{i=1}^s\m_i,\qquad \forall |S|=t.
\eeq
If we remove every edge from a maximal collection of edge disjoint
copies of $H_i,\,i=1,2,\ldots,\ell$ then we destroy all copies of $H_i,i=1,2,\ldots,  \ell$.
Furthermore, no $t$-set will be independent if
$$E_0>10\sum_{i=1}^s\m_i.$$
This is equivalent to
$$10\sum_{i=1}^s\binom{t}{r}\binom{n}{v_i-r}e_ip^{e_i}\frac{r!(v_i-r)!}{\a_i}<
\frac12\binom{t}{r}p$$
or
$$\sum_{i=1}^s\binom{n}{v_i-r}e_ip^{e_i-1}\frac{r!(v_i-r)!}{\a_i}<\frac{1}{20}$$
and this is implied by
$$\sum_{i=1}^s\frac{e_ic_1^{e_i-1}}{\a_i}<\frac{1}{20r!}.$$
This follows from the definition of $c_1$ and so
\beq{8}
after\ removal\ of\ edges,\ \a(G)\leq t.
\eeq

Thus the chromatic number is at least
$$k=\frac{n}{t}.$$
We re-express things to eliminate $n$. We have
$$k=\frac{1}{c_2}\bfrac{c_1}{r!\log
  n}^{1/(r-1)}n^{1-1/(\r(r-1))},$$
$$k^{(r-1)/(r-1-{1/\r})}=\frac{1}{c_2^{(r-1)/(r-1-{1/\r})}}\bfrac{c_1}{r!\log
  n}^{1/(r-1-{1/\r})}n.$$
Now we see from this that
$$\frac{r-1}{r-1-{1/\r}}\log k\sim \log n.$$
(Here $\sim$ denotes $=(1+o(1))$ as $k,n\to\infty$.)

So,
$$n\sim k^{(r-1)/(r-1-{1/\r})}c_2^{(r-1)/(r-1-{1/\r})}\bfrac{r!(r-1)\log
  k}{c_1(r-1-{1/\r})}^{1/(r-1-{1/\r})}.$$
Substituting in \eqref{1}, this gives
\begin{multline*}
E_p\leq \frac{2c_1}{r!}\brac{k^{(r-1)/(r-1-{1/\r})}c_2^{(r-1)/(r-1-{1/\r})}\bfrac{r!(r-1)\log
  k}{c_1(r-1-{1/\r})}^{1/(r-1-{1/\r})}}^{r-1/\r}\\
=2\bfrac{r!}{c_1}^{1/(r-1-{1/\r})}c_2^{(r-1)(r-
{1/\r})/(r-1-{1/\r}))}\bfrac{r-1}{r-1-{1/\r}}^{(r-{1/\r})/(r-1-{1/\r})}\\
\times k^{(r-1)(r-{1/\r})/(r-1-{1/\r})}(\log
k)^{(r-{1/\r})/(r-1-{1/\r})}.
\end{multline*}

Note that $n\to\infty$ implies $k\to\infty$ and so this completes the proof of Theorem \ref{main}.
\proofend
\section{$(r,\ell)$-systems}
In this section we give the short poof of Theorem \ref{rl}.  The only observation we need, which is very simple, is that an $(r,\ell)$-system is
one where a particular finite list of hypergraphs is forbidden.
\medskip

\noindent
{\bf Proof of Theorem \ref{rl}.}  We use Theorem \ref{main}.   Let
$H_i,\,i=1,2,\ldots,r-\ell$ be the hypergraph consisting of two edges intersecting in
$\ell+i-1$ vertices.  Then an $(r,\ell)$-system is one which is $\cH$-free, where $\cH=\{H_1, \ldots, H_{r-\ell}\}$. Using the notation of the previous section we have
\begin{eqnarray*}
\r&=&\frac{1}{r-\ell}\\
s&=&1\\
\a_1&=&\ell!(r-\ell)!^2\\
c_1&=&\frac{\ell!(r-\ell)!^2}{100r!}\\
c_2&=&10^{1/(r-1)}.
\end{eqnarray*}
Plugging these values into the expression for $c_\cH$ in Theorem \ref{main} gives us our expression for $c_{r,l}$.
This completes the proof of Theorem \ref{rl}.
\proofend
\section{Independent neighborhoods}

In this section we prove Theorem \ref{indnbd}.  We need the following three Lemmas.
The first was proved in \cite{EL}
and follows immediately from the Local Lemma.

\begin{lemma} {\bf (\cite{EL})}\label{LL}
Let $r \ge 2$ and let $G$ be an $r$-graph with maximum degree at most $k^{r-1}/4r$.  Then the chromatic number of $G$ is at most $k$.
\end{lemma}

The next Lemma has been proved by several researchers. In the form below it essentially appears in \cite{KMRT}.

\begin{lemma} \label{indchrom}
Let $0<\alpha\le 1/2$ and let $G$ be a hypergraph on $n$ vertices.  Suppose that every subhypergraph $P$ of $G$ (including $G$ itself) with $m$ vertices has an independent set of size $m^{\alpha}$.  Then $G$ has chromatic number at most $2n^{1-\alpha}$.
\end{lemma}

Our final Lemma is fairly straightforward, and generalizes the easy argument that an $n$ vertex triangle-free graph has an independent set of size at least $\sqrt n$ (actually, much more is guaranteed for graphs).

\begin{lemma} \label{indind}
Let $r \ge 3$ and let $G$ be an $n$-vertex $r$-graph with independent neighborhoods.  Then $G$ has an independent set of size at least $n^{1/r}$.
\end{lemma}

\noindent
{\bf Proof.}  Let $\Delta$ be the maximum size of a neighborhood of an $(r-1)$-set of vertices. Then
$$\sum_{v \in V(G)}d(v)=\sum_{|S|=r-1} d(S)\le {n \choose r-1}\Delta<\frac{\Delta n^{r-1}}{2}.$$
Consequently, the average degree $d$ of $G$ satisfies  $d \le \Delta n^{r-2}/2$.  Now by Tur\'an's theorem, $G$ has an independent set of size at least $(1-1/r)n/d^{1/{(r-1)}}$.  Therefore, we have an independent set of size at least $\max\{\Delta, (1-1/r)n/d^{1/{(r-1)}}\}\ge n^{1/r}$.
\qed

\noindent
{\bf Proof of Theorem \ref{indnbd}.}
For the upper bound, we apply Theorem \ref{main} with $\cH=\{F_r\}$.  In the notation of the proof of Theorem \ref{main}, we have
\begin{eqnarray*}
\r&=&2\\
s&=&1\\
\a_1&=&(r-1)!^2\\
c_1&=&\bfrac{(r-1)!^2}{50r!(r+1)}^{1/r}\\
c_2&=&10^{1/(r-1)}.
\end{eqnarray*}
Plugging these values into the expression for $c_\cH$ in Theorem \ref{main} gives us our expression for $c_\cI$.
For the lower bound, suppose that $G$ is an $r$-graph with independent neighborhoods and $|G|=bk^{r+1/(r-1)}$
where $b=1/(40r^22^r)$. Let $k$ be sufficiently large and even (a similar argument works for odd $k$) and let $A$ be the set of vertices in $G$ with degree
less than $d=k^{r-1}/(2r2^r)$.   By Lemma \ref{LL}, we can color the induced subhypergraph $G[A]$ properly by $k/2$ colors.   Let $G'\subset G$  be the $r$-graph induced by the uncolored vertices.  Since every vertex of $G'$ has degree (in $G$) at least $d$, the number of vertices $n$ of $G'$ satisfies $n\le rbk^{r+1/(r-1)}/d<k^{r/(r-1)}/20$.  Applying Lemmas \ref{indind} and \ref{indchrom}, we conclude that $G'$ has a  proper coloring where the number of colors is at most
$$2n^{1-1/r}<(2/5)k< k/2.$$
Putting these two colorings together yields a proper coloring of $G$ with at most $k$ colors.
\qed

\section{Excluding a hypertree}
In this section we prove Theorem \ref{Tmain}.  Recall that an
$r$-tree is a connected $r$-forest, where connected means that for
every two vertices $x,y$, there is a sequence of edges $E_1, \ldots,
E_l$ such that $x \in E_1$, $y \in E_l$, and $E_i \cap E_{i+1}\ne
\emptyset$ for all $i=1,2, \ldots, l-1$. If \( T \) is an $r$-tree,
then an edge \(e \in T \) is a {\em leaf} if \( e \) contains at
most one vertex of degree greater than one.
\bigskip

\noindent {\bf Proof of Theorem \ref{Tmain}.} We begin by
inductively defining a sequence of collections of $r$-trees.  Set \(
\cF_0 = \{ T \} \).  For \( i =1, \dots ,t-1 \) let \( \cF_i \) be
the collection of $r$-trees given by deleting a leaf from some
$r$-tree in  \( \cF_{i-1} \). Given an $r$-tree $T \in \cF_i$ with
$i \ge 1$, say that a vertex $v$ of $T$ is a {\em connector} if
adding a leaf  to $v$ results in a tree $T' \in \cF_{i-1}$.  Such a
connector exists by the way we have defined the sequence
$\{\cF_i\}$. Note that each $r$-tree in \( \cF_{i} \) has \( t - i
\) edges and spans \( (r-1)(t-i) + 1\) vertices.

Let \(V\) be the vertex set of the $r$-graph \(G\) that does not
contain a copy of \( T \). We use the sets \( \cF_1, \dots ,
\cF_{t-1} \) to define a collection of disjoint subsets of \( V \).
Set \( G_1 = G \) and let \(A_1 \) be the set of vertices \(v \in
V\) with the property that there exists some $r$-tree \(T^\prime \in
\cF_1 \) such that \(G \) contains a copy of \( T^\prime \) with \(
v \) as a connector. For each such vertex \(v\) let \( X_v \) be the
set of vertices (other than \(v\)) spanned by {\bf one} of these
$r$-trees \( T^\prime \in \cF_1 \) that contain \(v\) as a
connector. Note that \( G_2 := G \left[ V \setminus A_1 \right] \)
does not contain any copies of any $r$-tree in \( \cF_1 \) (such a
copy would include a connector and all such vertices were gathered
into \(A_1\)).

Now, suppose disjoint sets \( A_1, \dots, A_{i} \subseteq V \) have been defined with the
following properties:
\begin{enumerate}
\item[(i)] If \( v \in A_j \) then there is a copy of an $r$-tree \( T^\prime \in \cF_{j} \) in \(G_{j}\) with \(v\)
as a connector.  The set of vertices (other than \(v\)) spanned by
{\bf one} such $r$-tree is \( X_v \).
\item[(ii)] The graph
\[ G_{i+1} = G \left[ V \setminus \cup_{j=1}^{i} A_j \right] \]
does not contain any copy of an $r$-tree in \( \cF_{i} \).
\end{enumerate}
Let \(A_{i+1}\) be the set of connectors of copies of $r$-trees in
\( \cF_{i+1} \) in \( G_{i+1} \).  For each \( v \in A_{i+1} \) let
\( X_v  \cup \{v\}\) be the vertex set of {\bf one} of the $r$-trees
\( T^\prime \in \cF_{i+1} \) that lies in \(  G _{i+1} \) and
contains \(v\) as a connector.  Note that $V-\cup_{i=1}^{t-1} A_i$
contains no edges of $G$, since $\cF_{t-1}$ is the tree with one
edge.

Now consider the graph \( H_G \) with vertex set \(V\) and edge set
\[ \bigcup_{v \in V} \left\{ \{ u,v \} : u \in X_v \right\}. \]
In words we put an edge between each vertex \(v\) and every vertex in the set
\( X_v \).  Note that every induced subgraph of \( H_G \) has average degree
bounded above \( 2(r-1)(t-1) \) (as all edges in the subgraph induced by \(Y\)
are `generated' by one of the vertices in \(Y\) and each such vertex `generates'
at most \( (t-1)(r-1) \) edges).  It follows that \( G_H\) is \( 2(r-1)(t-1) \)-degenerate
and can be colored with \( 2(r-1)(t-1) + 1 \) colors.  Let \( f \) be a proper coloring of \( H_G \)
with \(2(r-1)(t-1) + 1 \) colors.

We claim that \(f\) is also a proper coloring of \( G \).  Let \( e
\) be an edge in \(G\).  Then we have observed above that $e$ must
intersect some $A_j$. Let \(v \in e \cap A_i \) where \(i\) is the
smallest index such that \( e \cap A_i \neq \emptyset \). Consider
the $r$-tree \( T^\prime \in \cF_i \) that contains \( v \) as a
connector and spans \( \{v\} \cup X_v \). If \( e \cap X_v =
\emptyset \) then $X_v \cup e$ spans a copy of an $r$-tree in \(
\cF_{i-1} \), which contradicts the properties of the sets \( A_1,
\dots, A_{t-1} \).  Therefore \( e \) contains some vertex \(u \in
X_v \). As \( \{ u,v\} \in E(H_G) \), \( f \) assigns \( u \) and
\(v\) different colors.\qed

\section{Cliques}
In this section we prove Theorem \ref{cliques}. We must provide a construction that has fewer edges than the one in Theorem \ref{main} when $r \ge 3$. It is motivated by similar constructions in Ramsey-Tur\'an theory.
\medskip

\noindent
{\bf Construction.}  Fix $r \ge 3$.  Let $G$ be the $r$-graph with vertex set $V=[n]$ obtained by the following random process.  For each $i \in [n]$, randomly partition $\{i+1, \ldots, n\}$ into $r-1$ sets $V^i_1, \ldots, V^i_{r-1}$, each of size
$\lceil\frac{n-i}{ r-1}\rceil$ or $\lfloor \frac{n-i}{r-1}\rfloor$.  Now add all edges of the form $\{i, v_1, \ldots, v_{r-1}\}$, where $i<v_j$ and $v_j \in V^i_j$ for all $j$.  \qed

\medskip

\noindent
{\bf Proof of Theorem \ref{cliques}.}  Let us first observe that $G$ contains no copy of $K_{r+1}^r$.
Indeed, if $K$  is such a copy, let $i$  denote its smallest vertex. Since there are $r$ other vertices in $K$, by the pigeonhole principle, two of these, say $w$ and $y$ lie in $V^i_j$ for some $j$.  But this means that there is no edge of $G$ containing all three of $i,w,y$, and in particular, at least one (in fact many) edge of $K$ is missing in $G$. This contradiction implies that $G$ contains no $(r+1)$-clique.

If $|G|$ denotes the number of edges in $G$, then by counting edges from their leftmost endpoint we see that
$$|G| \le \sum_{i=1}^{n-r+1}\left(\frac{n-i}{r-1}+1\right)^{r-1}
<\frac{1}{(r-1)^{r-1}}\sum_{j=1}^{n}j^{r-1}<\frac{n^r}{(r-1)^{r-1}}.$$
Let us obtain an  upper bound on the independence number of $G$. For any $r$-tuple $f=\{i_1,\ldots,i_r\}$ with $i_1<\cdots<i_r$, let $\cE_f$ be the event that $f \in G$.  If $f=\{i_1,\ldots,i_r\}$ and $f'=\{i_1',\ldots,i_r'\}$ with $i_1<i_2<\cdots < i_r$ and $i'_1<i'_2<\cdots i'_r$,
then
\medskip
\beq{star}
\cE_f\ and\ \cE_{f'}\ are\ independent\ if\ and\ only\ if\ i_1\ne i_1'.
\eeq

\noindent
Now pick a set $S=\{v_1, \ldots, v_s\} \subset V$ with $v_1<v_2 <\cdots <v_s$.  Let $G_i$ be the set of edges in $G[S]$ whose smallest vertex is $v_i$.  Then
\begin{eqnarray}
\Pr(S \hbox{ is independent})&=&\Pr(G_i=\emptyset \hbox{ for all $i=1, \ldots, s$}) \nonumber\\
&<&\prod_{i=1}^{s/2}\Pr(G_i=\emptyset)\nonumber\\
&=&\prod_{i=1}^{s/2}\Pr(\exists j:\;V^i_j\cap \set{v_{i+1},\ldots,v_s}=\emptyset)\nonumber\\
&<&\prod_{i=1}^{s/2}\frac{(r-1)\binom{(r-2)\lceil(n-v_i)/(r-1) \rceil}{s-i} }{\binom{n-v_i}{s-i}}\nonumber\\
&<& \prod_{i=1}^{s/2} r \left(\frac{r-2}{r-1}\right)^{s-i}\nonumber\\
&<&\left(re^{-s/(2r)}\right)^{s/2} \nonumber
\end{eqnarray}
where the first inequality holds due to \eqref{star}.
Consequently, the expected number of independent sets of size $s$ in $G$ is at most
$$\binom{n}{s}\cdot  \left(re^{-s/(2r)}\right)^{s/2}<\brac{\frac{ner^{1/2}}{s}e^{-s/(4r)}}^s<1$$
as long as $s>4r\log n$.   This shows that there exists such a $G$ with chromatic number $k$ at least  $n/(4r\log n)$.  Since the number of edges in $G$ is at most $\frac{n^r}{(r-1)^{r-1}}$, this construction gives
$$m_k(K_{r+1}^r)<d_r(k\log k)^r$$
where $d_r\leq 5^rr^r/(r-1)^{r-1}$.

Now we prove the statement about graphs.  By standard results in Ramsey theory, every $K_s$-free graph on $n$ vertices has an independent set of size at least $n^{\delta_s}$, where $\delta_s>0$. Now choose $\e=\e_s$ such that $0<\e<1/(1-\delta_s)-1$.  Suppose that
$G$ is a $K_s$-free graph with independent neighborhoods and $k^{2+\e}$ edges, where $k$ is sufficiently large.  Let $A$ be the set of vertices in $G$ with degree
less than $k/2-1$.   We can greedily color the induced subgraph $G[A]$ properly by $k/2$ colors.   Let $G'\subset G$  be the subgraph induced by the uncolored vertices.  Since every vertex of $G'$ has degree (in $G$) at least $k/2$, the number of vertices $n$ of $G'$ satisfies $n\le 4k^{1+\e}$.  By the choice of $\delta_s$, every $m$-vertex subgraph of $G'$ has an independent set of size at least $m^{\delta_s}$.
Hence by Lemma~\ref{indchrom}, we conclude that $G'$ has a  proper coloring where the number of colors is at most
$$2n^{1-\delta_s}<2(4k^{1+\e})^{1-\delta_s}<k/2,$$
where the last inequality holds by the choice of $\e$ and the fact that $k$ is sufficiently large.
Putting these two colorings together yields a proper coloring of $G$ with at most $k$ colors.\qed
\section{Concluding remarks and open problems}
In this section we repeat some of the open questions  mentioned throughout the paper and state a couple of new ones as well.
\bigskip

\noindent
$\bullet$ Attempts to improve the lower bound in Theorem \ref{indnbd} lead to the following question which is independently interesting. Suppose that $G$ is an $r$-graph with independent neighborhoods and maximum degree $\Delta$. What are the best upper bounds one can obtain on the chromatic number of $G$?   The Local lemma gives $O({\Delta}^{1/(r-1)})$, but the results for
 the graph case ($r=2$) suggest that one should be able to improve this to
$O(({\Delta}/\log {\Delta})^{1/(r-1)})$. The $r=2$ case, that
triangle-free graphs with maximum degree $\Delta$ have chromatic
number at most $O(\Delta/\log{\Delta})$, is a deep result due to
Johansson, but those ideas do not extend to $r>2$. When $r=3$ we
pose the following weaker statement.

\medskip
\noindent{\bf Problem.}  Let $G$ be a 3-graph with independent neighborhoods and maximum degree $\Delta$.  Prove that the chromatic number of $G$ is $o(\sqrt{\Delta})$.
\medskip

\noindent
A much stronger statement for graphs has been conjectured  by Alon-Krivelevich-Sudakov~\cite{AKS}.

 As we mentioned earlier, we do not believe that the order of magnitude of the upper bound in Theorem \ref{indnbd} is correct either. Perhaps some generalization of Kim's construction for $R(3,t)$ would improve the log factors.

\bigskip

\noindent $\bullet$ Let $T$ be an $r$-tree with $t$ edges and $G$ be
an $r$-graph containing no copy of $T$.  When $r=2$, it is
well-known that the chromatic number of $G$ is at most $t$, and this
is sharp.  Theorem \ref{Tmain} gives an upper bound of about $2rt$,
but again the best lower bound we have is roughly $t$.  It would be
very interesting to narrow this gap, in particular to determine
whether the coefficient of $t$ depends on $r$ in an essential
way\footnote{Recently Po-Shen Loh has proved optimal results for
this problem}.

\bigskip

\noindent
$\bullet$  Our final question is perhaps too ambitious given the current state of knowledge, and pertains to Theorem \ref{cliques}.

\medskip
\noindent{\bf Problem.}
Characterize all 3-graphs $H$ such that $m_k(H)=k^{3+o(1)}$.

\end{document}